\documentclass[oneside]{article}
\usepackage{amsfonts} % for ... ?
\usepackage{listings}
\newcommand{\frc}{\ensuremath{\mbox{{\it frac}}}}

{}
{}
{}
{}
{}
{}

\title{{\bf Egyptian Multiplication\\
            and some of its ramifications}}
\author{M.H. van Emden\\
        {\small Technical Report DCS-362-IR}\\[-1mm]
        {\small Department of Computer Science}\\[-1mm]
        {\small University of Victoria}
       }
\date{}

\begin{document}
\maketitle
\date{}

\begin{abstract}

Multiplication and exponentiation can be defined by equations
in which one of the operands is written as the sum of powers
of two. When these powers are non-negative integers, the
operand is integer; without this restriction it is a fraction.
The defining equation can be used in {\sl evaluation mode}
or in {\sl solving mode}.
In the former case we obtain ``Egyptian'' multiplication,
dating from the 17th century BC.
In solving mode we obtain an efficient algorithm for division
by repeated subtraction, dating from the 20th century AD.
In the exponentiation case we also distinguish between
evaluation mode and solving mode.
In the former case we obtain an algorithm for fractional
powers;
in the latter case we obtain an algorithm for logarithms,
the one invented by Henry Briggs in the 17th century AD.
\end{abstract}

\section{{\large Egyptian multiplication}}

To begin at the beginning,
let the starting point be the multiplication algorithm
in the Rhind papyrus, a document found in Egypt.
My source is an article by James R. Newman \cite{newmanRhind}.
The document in question was a scroll of 18 feet long
written in about 1700 B.C.
Though not found intact, enough of the parts have been
compiled to obtain a coherent whole.
The scroll is a collection mathematical exercises and
practical examples. Most of it is of no mathematical
interest because of the unwieldy notation for fractions.
However, one algorithm, for multiplying two integers,
stands out for its clarity, simplicity,
and relevance for computer programming.

Newman illustrates the algorithm by multiplying 23 and 27 as
shown in Figure~\ref{prog:rhind}.
\begin{figure}[h]
\hrule \vspace{2mm}
\begin{center}
\begin{minipage}[t]{4in}
\begin{verbatim}

                 \   1     27
                 \   2     54
                 \   4    108
                     8    216
                 \  16    432
              ---------------
            Total   23    621

\end{verbatim}
\end{minipage}
\end{center}
\caption{
\label{prog:rhind}
Multiplication as found in the Rhind papyrus.
The total includes only the rows marked by a backslash.
These rows correspond to the digit 1 in {\tt 10111}, the
binary representation of 23.
}
\vspace{2mm}
\hrule
\end{figure}
It is not clear from his article whether this example actually
occurs in the papyrus or whether it is  his way of presenting
the idea gleaned from the papyrus.
Let us translate the idea to modern mathematical
notation.

The idea is that $a\times b$,
with non-negative integer $b$, can be written as
$a\times \sum_{i=0}^\infty d_i 2^i$
where the $d_i$ are $0$ or $1$.
Of course, from a certain $i$ onward, all the $d_i$ are zero.
We have
$a\times b = 
 a\times \sum_{i=0}^\infty d_i 2^i =
 \sum_{i=0}^\infty d_i 2^i a
$.
In the second equality, the multiplication is converted
to additions (in the form of doublings).

In Newman's example $a=27$ and $b=23$, so that we get
$$
27\times 23 = 27+54+108+432
$$
corresponding to the binary expansion of $b$ being $10111$.
When we determine this expansion in the usual way,
the digits are generated in reverse order,
so that we get
$$
27\times 23 = 1\times 27+1\times 54+1\times 108+0\times 216
            + 1\times 432.
$$
This suggests evaluating $\sum_{i=0}^\infty d_i 2^i a$
by combining the generation of the binary digits with
the doublings of $a$ in a single loop, as in
\begin{verbatim}
    r := 0;
    while (b>0) {
      if (odd(b)) {r := r+a; b := b-1}
      a := a+a; b := b/2
    }
\end{verbatim}
with the invariant $a_0\times b_0 = r+a\times b$
where 
$a_0$ and $b_0$ are the initial values of $a$ and $b$.

Although Newman's presentation suggests that a table of
the successive doublings be kept, the above code shows
that each result of doubling can overwrite the previous
one.

This is an example of using the binary expansion to determine
the value of $a\times b$ when both operands are known.
Let us call this ``evaluation mode''.
We can also use the binary expansion in ``solving mode'':
to determine an unknown $b$ when $a$ and $c$ are given in
$a\times b=c$.
As this does not in general have a solution in integers,
we solve instead $a\times b=c-r$ with $0\leq r < b$.
In solving mode we compute the inverse, which in this case
is division:
$$c-r = a\times b = 
 a\times \sum_{i=0}^\infty d_i 2^i =
 \sum_{i=0}^\infty d_i 2^i a
$$.

For example
$$
626-r = d_0 27+d_1 54+d_2 108+d_3 216+d_4 432.
$$
Here it is apparent that $d_4=1$,
so that we can expose $d_3$ by subtracting $432$ from both
sides.
$$
194-r = d_0 27+d_1 54+d_2 108+d_3 216
$$
from which we conclude the $d_3=0$ and $d_2=1$.
Further subtractions give
$
86-r = d_0 27+d_1 54
$
and
$
32-r = d_0 27
$
from which  we conclude that $d_0=1$ and $r=5$,
so that division of $626$ by $27$ gives quotient
$11101$ (in binary, with most significant digit first)
and leaves a remainder of 5.

\paragraph{}
The presentation by Newman suggests that a table of doublings
of $a$ be kept. This is not necessary, as the next doubling
can be written over the current one. 
Determining the binary expansion of $b$ tells one when to stop.

Most of the Rhind papyrus is concerned with division.
No similarly simple and clean algorithm is found by Newman.
This is probably also the case for the study \cite{chace}
by Arnold Chace, on which Newman's article is based.

The above numerical example for integer division
shows that one needs to know
in advance the largest doubling needed. Once this is known,
the lesser doublings are efficiently obtained by halving;
again no table is needed.
See Figure~\ref{prog:deBruijn}.

\begin{figure}[h]
\begin{center}
\begin{minipage}[t]{4in}
\hrule \vspace{2mm}
\begin{verbatim}
void qr1(int a, int d, int* Q, int* R) {
  int r, dd, q;
  r = a; dd = d; q = 0; // a = q*dd + r
  while (dd <= r) dd = 2*dd;
  // dd = 2^i * d with least i s.t. dd>r & a = q*dd + r
  while (dd != d) {
    // a = q*dd + r
    dd = dd/2; q = 2*q;
    // a = q*dd + r
    if (dd <= r) { r = r-dd; ++q; }
    // a = q*dd + r
  }
  *Q = q; *R = r; return;
}
\end{verbatim}
\hrule
\end{minipage}
\end{center}
\caption{
\label{prog:deBruijn}
Computing quotient and remainder.
Attributed in 1970 to N.G. de Bruijn by E.W. Dijkstra in
``Notes on Structured Programming'' EWD249, Section 5,
Remark 1. Reprinted in \cite{daDiHo}.
}
\end{figure}

\paragraph{}

This completes all that is said in this paper about multiplication
and  division. From now on the topic is exponentiation and its
inverse.
Just as multiplication can be done by repeated addition,
exponentiation can be done by repeated multiplication.
It can be speeded up in a similar way by using the binary expansion
of one of the operands. We will be concerned with the case of
a fractional operand, so that the binary expansion extends into
negative powers of two.
Again, algorithms can be obtained by using the defining
equation in the two modes used earlier: evaluating mode and 
solving mode.
We will need to take square roots because
what halving is to doubling, taking the square root is to
squaring.
There is an algorithm for square root that is of interest in its
own right.

\section{{\large The miracle of the square root}}
To obtain the square root $x$ of $a$, that is, to find $x$
such that $x^2 = a$, observe that this equality implies
$x = a/x$, hence $x = (x + a/x)/2$.
This suggests considering the sequence $x_0, x_1, \ldots$
defined by $x_{n+1} = (x_n + a/x_n)/2$ with some arbitrarily 
chosen $x_0$, say, 1. It is clear that if $x_n<\surd a$, then
$a/x_n > \surd a$, and vice versa.
If $x_n$ is regarded as a guess at the value of $\surd a$,
then $x_{n+1} = (x_n + a/x_n)/2$ is a plausible way of
getting a better guess.
It is not only plausible, but is guaranteed to converge
to $\surd a$. Moreover, the number of correct figures doubles
at every iteration.
When I speak of the {\sl miracle} of square root,
I have in mind this
combination of plausibility and algorithmic effectiveness.
\begin{figure}[h]
\begin{center}
\begin{minipage}[t]{4in}
\hrule \vspace{2mm}
\begin{verbatim}
double heron(double a) {
  double x = 1.0, newx, diff, eps = 1.0e-16;
  do { newx = (x + a/x)/2.0;
       diff = x-newx; x = newx;
  } while (-eps >= diff || eps <= diff);
  return newx;
}
\end{verbatim}
\hrule
\end{minipage}
\end{center}
\caption{
\label{prog:heron}
Heron's algorithm as function in C.
The keyword {\tt double} denotes the type of double-length
floating-point number.
}
\end{figure}

The method is an instance of Newton's method of finding
roots of non-linear equations, and is often presented as such.
This is a pity because it suggests that one needs calculus
to understand a good method for the square root.
This is not the case, as the algorithm is described by Heron
of Alexandria who lived two thousand years ago\footnote{
There is ``informed conjecture'' that this algorithm
was used by the Babylonians, and thus may be roughly as old as
Egyptian multiplication.
}.

\section{{\large Fractional powers}}
The Heron algorithm might suggest looking for an equally
elegant and effective algorithm for
cube roots, fifth roots, $\ldots$.
In the context of this paper a more attractive option
is to exploit square roots
for raising $a$ to any power
of the form $p/q$, for non-negative integer $p$ and any
positive integer $q$.

The mathematical basis of the algorithm is
\begin{eqnarray*}
a^{2p/q} &=& (a*a)^{p/q} \\
a^{p/q} &=& a(a^{(p-q)/q}) \\
a^{p/q} &=& (a*a)^{p/2q} \\
a^{p/q} &=& (\surd a)^{2p/q} \\
\end{eqnarray*}
These equalities suggest mutual adjustments among $a$, $p$,
and $q$ in such a way that after a sufficient number of steps
the desired power is trivial to obtain.
See Figure \ref{prog:fractPowSP}.
\begin{figure}[h!]
\begin{center}
\begin{minipage}[t]{4in}
\hrule \vspace{2mm}
\begin{verbatim}
double fp(double a, int p, int q) {
// fractional power: returns a^(p/q)
  double z, eps = 1.0e-14;
   z = 1.0;
  while (p>q) {
    if ((p%2) == 0) { p = p/2; a = a*a; }
    if ((p%2) == 1) { p = p-q; z = z*a; }
  }
  while (1) {
    if (1.0-eps < a && a < 1.0+eps) return(z);
    if (p == q) return(z*a);
    if (p < q) { a = heron(a); p = 2*p; continue; }
    if (p > q) { p = p-q; z = z*a; continue; }
  }
}
\end{verbatim}
\hrule
\end{minipage}
\end{center}
\caption{
\label{prog:fractPowSP}
Algorithm to compute power with exponent given as
ratio of positive integers.
}
\end{figure}

Fractional powers are defined by $a^{p/q}=b$.
To compute such a power, 
$a$, $p$, and $q$ are given while $b$ is unknown.
This defining equation is thus used in evaluation mode.
Ideally, the same defining equation can be used in solving
mode, with $a$ and $b$ given.
But an unknown in the form $p/q$ is awkward.
Therefore we use $a^t=b$ as defining equation,
where $t$ is understood to be a fraction.

In Figure \ref{prog:fractPowFLPT} we list a C function
to compute $a^t$
even though this function is merely a predictable variant of
the one listed in Figure \ref{prog:fractPowSP}.

\begin{figure}[h!]
\begin{center}
\begin{minipage}[t]{4in}
\hrule \vspace{2mm}
\begin{verbatim}
double fpFLPT(double a, double t) {
// fractional power: returns a^t
  double z, eps = 1.0e-14;
  z = 1.0;
  // a = a0 & t = t0
  // maintain a0^t0 = z*a^t
  while (t > 1.0) { t = t/2; a = a*a; }
  // t <= 1.0   t will dance around 1
  while (a < 1.0-eps || 1.0+eps < a) {
    if (t >= 1.0) { t = t-1.0; z = z*a; }
    if (t <  1.0) { t = 2*t; a = heron(a); }
  }
  return(z);
}
\end{verbatim}
\hrule
\end{minipage}
\end{center}
\caption{
\label{prog:fractPowFLPT}
Algorithm to compute power with exponent given as
floating-point number.
}
\end{figure}

\section{{\large In the footsteps of Henry Briggs}}

The Rhind papyrus uses a table of doublings of the
multiplier. A modern version of this algorithm
creates the table implicitly, each entry overwriting
the previous one.
In this section we will see that a similar table is
useful in computing logarithms.
This use occurred early in the history of
computation; barely a quarter of a century after Simon Stevin
taught the world to use decimal positional notation for fractions.
I quote from {\sl The Feynman Lectures on Physics} (\cite{feynmLect},
page 22-6):
\begin{quote}
{\sl
This is how logarithms were originally computed by Mr Briggs of
Halifax, in 1620. He said ``I computed successively 54 square
roots of 10.''
}
\end{quote}
Briggs might easily have picked another number of square
roots.
On my computer (and on yours probably as well) the iteration
\begin{verbatim} 
      double b = sqrt(10.0); while (b>1.0) b = sqrt(b);
\end{verbatim}
stops after producing 53 distinct square roots.
I don't know how Briggs decided on 54, but this
occurrence of ``53'' has an explanation. The {\tt double}
type in C is represented by the double-length format of
the IEEE floating-point standard, which has a mantissa of
52 bits. This might suggest a relative precision of 52 bits,
were it not for the fact that
the first bit of the significand is suppressed in the format
because it is always 1.
Because that is taken into account,
one gets effectively a precision of 53 bits.
A curious near-coincidence with the choice of Briggs in 1620.
For future reference we note that
53 bits is equivalent to about 16 decimal places.

We let the Feynman lectures (\cite{feynmLect}, page 22-7) continue
with
\begin{quote}
{\sl
$\ldots$ he calculated sixteen decimal places,
and then reduced his answer to fourteen places when he
published it, so that there were no rounding errors.
He made tables of logarithms to fourteen decimal places by
this method, which is quite tedious.
But all logarithm tables for three hundred years were borrowed
from Mr Briggs's tables by reducing the number of decimal places.
Only in modern times, with the WPA\footnote{
Works Progress Agency, a US government agency active in the
1930s to carry out public works,
thus alleviating unemployment.
}
and computing machines, have
new tables been independently computed.
}
\end{quote}

The role of tables so far has been conceptual:
no data structure was created to be filled with results
of squarings or doublings, of square roots or 
halvings.
Instead, the contents of the notional tables
were constructed on the fly and overwritten
promptly after use.

Rather than pursue an algorithm for base-10 logarithm,
we note that Briggs's method of computing a table of
iterated square roots works equally well for any base $b$
(within a suitable range) as it does for base 10.

\section{{\large Logarithms to any base}}

So far, when we considered
fractional powers, as in $b^x = a$ with fractional $x$,
we only used this equality in evaluation mode.
In solving mode, $b$ and $a$ are given and one is to find $x$.
In other words, one is to find the logarithm of $a$ to base $b$.
We determine
the $x$ in $b^x = a$ when $a$ and $b$ are given,
with $b>1$ and $1 \leq a \leq b$.
From these assumptions we conclude that $0 \leq x \leq 1$.
Accordingly, we have that
$$
x = \sum_{i=1}^{53} d_i 2^{-i},
$$
with $d_i \in \{0,1\}$
and uniquely determined in the precision anticipated for $x$.

A plausible way to attain the goal of $b^x = a$ is to suppose
we have a $k$ such that
$$
x = \sum_{i=1}^{k-1} d_i 2^{-i} + d_k 2^{-k} + \sum_{i=k+1}^{53} d_i 2^{-i}.
$$
Initially this is easy to make true with $k=1$.
To increase $k$ by one
we maintain program variables $z$ and $\frc$ such that
$$\log_b z = \sum_{i=1}^{k-1} d_i 2^{-i} \mbox{ and } \frc = 2^{-k}.  $$
To discover whether $d_k = 1$ or $d_k = 0$ 
we test $z*b < a$. Truth indicates the former, falsity the latter.
With $b=1$, $x$ cannot increase any more,
and we know that its true value (to the anticipated precision)
is representable as a sum of powers
of 2, and we have included all of the powers that should be included.

Such considerations lead to the function listed in
Figure~\ref{prog:bLog}.

\begin{figure}[h!]
\begin{center}
\begin{minipage}[t]{4in}
\hrule \vspace{2mm}
\begin{verbatim}
double lg(double b, double a) {
// Precondition: 1 <= a < b
// Returns the logarithm of a to the base b.
  double z = 1, frac = 1, x = 0;
  while (b > 1.0) {
    b = heron(b); frac /= 2;
    if (z*b < a) { z *= b; x += frac; }
  }
  return x;
}
\end{verbatim}
\hrule
\end{minipage}
\end{center}
\caption{
\label{prog:bLog}
A function for computing logarithms to base $b$.
}
\end{figure}

%\appendix

\section*{{\large Acknowledgements}}
Thanks to Paul McJones for corrections and discussions.

\bibliographystyle{plain}
\bibliography{bibfile}

\end{document}